\documentclass[12pt]{amsart}

\usepackage[francais,english]{babel}
\usepackage[T1]{fontenc}
\usepackage[latin1]{inputenc}
\usepackage{amsmath,amssymb,amsthm,mathrsfs}
\usepackage[all]{xy}
\usepackage{graphicx}
\usepackage{enumerate}

\numberwithin{equation}{section}

\DeclareMathOperator{\Fct}{Fct}
\DeclareMathOperator{\smd}{smd}
\DeclareMathOperator{\add}{add}
\DeclareMathOperator{\RHom}{RHom}
\DeclareMathOperator{\Ext}{Ext}
\DeclareMathOperator{\fRHom}{R\mathcal{H}om}

\DeclareMathOperator{\Hom}{Hom}

\DeclareMathOperator{\id}{id}
\DeclareMathOperator{\opp}{op}
\DeclareMathOperator{\Mod}{Mod}
\DeclareMathOperator{\sMod}{\mathfrak{M}od}

\DeclareMathOperator{\Hn}{H}

\DeclareMathOperator{\fHom}{\mathcal{H}om}

\DeclareMathOperator{\gr}{gr}
\DeclareMathOperator{\Rg}{R\Gamma}
\DeclareMathOperator{\dR}{R}

\begin{document}

\newtheorem{thm}{Theorem}[subsection]
\newtheorem*{theorm}{Theorem}
\newtheorem{conjecture}[thm]{Conjecture}
\newtheorem*{Notation}{Notation}
\newtheorem{defi}[thm]{Definition}
\newtheorem{prop}[thm]{Proposition}
\newtheorem{cor}[thm]{Corollary}
\newtheorem{lemme}[thm]{Lemma}
\newtheorem{Rem}[thm]{Remark}
\newtheorem{Ex}[thm]{Exemple}

\newcommand{\On}[1]{\mathcal{O}_{#1}}
\newcommand{\En}[1]{\mathcal{E}_{#1}}
\newcommand{\Fn}[1]{\mathcal{F}_{#1}} 
\newcommand{\tFn}[1]{\mathcal{\tilde{F}}_{#1}}
\newcommand{\hum}[1]{hom_{\mathcal{A}}({#1})}
\newcommand{\hcl}[2]{#1_0 \lbrack #1_1|#1_2|\ldots|#1_{#2} \rbrack}
\newcommand{\hclp}[3]{#1_0 \lbrack #1_1|#1_2|\ldots|#3|\ldots|#1_{#2} \rbrack}
\newcommand{\catMod}{\mathsf{Mod}}
\newcommand{\Der}{\mathsf{D}}
\newcommand{\Ds}{D_{\mathbb{C}}}
\newcommand{\DG}{\mathsf{D}^{b}_{dg,\mathbb{R}-\mathsf{C}}(\mathbb{C}_X)}
\newcommand{\lI}{[\mspace{-1.5 mu} [}
\newcommand{\rI}{] \mspace{-1.5 mu} ]}
\newcommand{\Ku}[2]{\mathfrak{K}_{#1,#2}}
\newcommand{\iKu}[2]{\mathfrak{K^{-1}}_{#1,#2}}
\newcommand{\Be}{B^{e}}
\newcommand{\op}[1]{#1^{\opp}}
\newcommand{\N}{\mathbb{N}}
\newcommand{\Ab}[1]{#1/\lbrack #1 , #1 \rbrack}
\newcommand{\Du}{\mathbb{D}}
\newcommand{\C}{\mathbb{C}}
\newcommand{\Z}{\mathbb{Z}}
\newcommand{\w}{\omega}
\newcommand{\K}{\mathcal{K}}
\newcommand{\Hoc}{\mathcal{H}\mathcal{H}}
\newcommand{\env}[1]{{\vphantom{#1}}^{e}{#1}}
\newcommand{\eA}{{}^eA}
\newcommand{\eB}{{}^eB}
\newcommand{\eC}{{}^eC}
\newcommand{\cA}{\mathcal{A}} 
\newcommand{\cB}{\mathcal{B}}
\newcommand{\cR}{\mathcal{R}}
\newcommand{\cL}{\mathcal{L}}
\newcommand{\cO}{\mathcal{O}}
\newcommand{\cM}{\mathcal{M}}
\newcommand{\cK}{\mathcal{K}}
\newcommand{\Hper}{\Hn^0_{\textrm{per}}}
\newcommand{\Dper}{\Der_{\textrm{per}}}
\newcommand{\Yo}{\textrm{Y}}
\newcommand{\gqcoh}{\textrm{gqcoh}}
\newcommand{\coh}{\textrm{coh}}
\newcommand{\cc}{\textrm{cc}}
\newcommand{\qcc}{\textrm{qcc}}
\newcommand{\qcoh}{\textrm{qcoh}}
\newcommand{\obplus}[1][i \in I]{\underset{#1}{\overline{\bigoplus}}}
\newcommand{\Lte}{\mathop{\otimes}\limits^{\rm L}}

\author{Fran\c{c}ois Petit}
\address{Institut de Math\'ematiques de Jussieu\\ UPMC\\ 4, place Jussieu\\ 75252 Paris Cedex 05, FRANCE.}
\email{fpetit@math.jussieu.fr}
\title{DG affinity of DQ-modules}
\begin{abstract}
In this paper, we prove the dg affinity of formal deformation algebroid stacks over complex smooth algebraic varieties. For that purpose, we introduce the triangulated category of formal deformation modules which are cohomologically complete and whose associated graded module is quasi-coherent. 
\end{abstract}
\maketitle
\section{Introduction}

Many classical results of complex algebraic or analytic geometry have a counterpart in the framework of Deformation Quantization modules (see \cite{KS3}). Let us mention a few of them, Serre duality, convolution of coherent kernels, the construction of Hochschild classes for coherent DQ-modules in \cite{KS3}, a GAGA type theorem in \cite{Gaga_DQ} and Fourrier-Mukai transforms in \cite{Mukai_Pantev} etc.

In this paper, we give a non-commutative analogous of a famous result of Bondal-Van den Bergh asserting the dg affinity of quasi-compact quasi-separate schemes (see \cite[Corollary 3.1.8]{BVdB}). In the framework of formal deformation algebroid stacks, the notion of quasi-coherent object is no more suited for this purpose. Thus, we introduce the notion  of cohomologically complete and graded quasi-coherent objects ($\qcc$ for short). The $\qcc$ objects of the derived category $\Der(\cA_X)$, where $\cA_X$ is a formal deformation algebroid stacks, form a full triangulated subcategory of $\Der(\cA_X)$ denoted $\Der_{\qcc}(\cA_X)$. This category can be thought as the deformation of $\Der_{\qcoh}(\cO_X)$ while deforming $\cO_X$  into $\cA_X$ (see Theorem \ref{deforq}). We prove that the image of a compact generator of $\Der_{\qcoh}(\cO_X)$ is a compact generator of $\Der_{\qcc}(\cA_X)$. The existence of a compact generator in $\Der_{\qcoh}(\cO_X)$ is granted by a result of Bondal-Van den Bergh (see loc. cit.). Hence, the category $\Der_{\qcc}(\cA_X)$ is dg affine.

The study of generators in derived categories of geometric origin has been initiated by Beilinson in \cite{Beil_gen}. The results of \cite{BVdB}  have been refined by Rouquier in \cite{Rou} where he introduced a notion of dimension for triangulated categories. Recently, in \cite{Toen_generateur} Toen generalized the results of Bondal and Van den Bergh and reinterpreted them in the framework of homotopical algebraic geometry.\\

This paper is organised as follows. In the first part, we recall some classical material concerning generators in triangulated category. We review, following \cite{KS3}, the notion of cohomological completeness and its link with the functor of $\hbar$-graduation. We finally state some results specific to deformation algebroid stacks on smooth algebraic varieties.

In the second part of the paper, we introduce the triangulated category of $\qcc$ objects, that is to say elements of $\Der(\cA_X)$ that are cohomologically complete and whose associated graded module is quasi-coherent. We prove that the category $\Der_{\qcc}(\cA_X)$ admits arbitrary coproducts. The coproduct is given by the cohomological completion of the usual direct sum (Proposition \ref{cocoqcc}) then we prove that $\Der_{\qcc}(\cA_X)$ is compactly generated (see Proposition \ref{genqcc} and Lemma \ref{impli}). Relying on a theorem of Ravenel and Neeman (see \cite{Rav} and \cite{Nee_comp})  we describe completely the compact objects of $\Der_{\qcc}(\cA_X)$ (see Theorem \ref{qccompact}). They are elements of 
$\Der^{b}_{\coh}(\cA_X)$ satisfying certain torsion conditions. Finally, we conclude this section by proving that $\Der_{\qcc}(\cA_X)$ is equivalent as a triangulated category to the derived category of a suitable dg algebra with bounded cohomology (see Theorem \ref{affinity}).

In the last section, we study $\qcc$ sheaves on an affine variety and prove that the equivalence of triangulated categories between $\Der^b_{\qcoh}(\cO_X)$ and $\Der^b(\cO_X(X))$ lifts to an equivalence between $\Der^b_{\qcc}(\cA_X)$ and the triangulated category $\Der^b_{cc}(\cA_X(X))$ of cohomologically complete $\cA_X(X)$-modules (see Theorem \ref{deforq}).

\noindent
\textbf{Acknowledgement.} I would like to express my gratitude to my advisor
P. Schapira for his patience and enthusiasm in sharing his mathematical knowledge. I would also like to thank Masaki Kashiwara for enlightening discussions and Grégory Ginot for his constant support.

\section{Review}

\subsection{Generators and Compactness in triangulated categories: a review}\label{envsec}

We start with some classical definitions. See \cite{Nee_book}, \cite{BVdB}.

\begin{defi}
Let $\mathcal{T}$ be a triangulated category. Let $\mathfrak{G}=(G_i)_{i \in I}$ be a set of objects of $\mathcal{T}$. One says that $\mathfrak{G}$ generates $\mathcal{T}$ if for every $F \in \mathcal{T}$ such that $\Hom_{\mathcal{T}}(G_i[n],F)=0$ for every $G_i \in \mathfrak{G}$ and $n \in \Z$ then $F \simeq 0$.
\end{defi}

Recall that if $\mathcal{T}$ is a triangulated category, then a triangulated subcategory $\mathcal{B}$ of $\mathcal{T}$ is called thick if it is closed under isomorphisms and direct summands.

\begin{defi}\label{epais}
Let $\mathcal{S}$ be a set of objects of $\mathcal{T}$. The smallest thick triangulated subcategory of $\mathcal{T}$ containing $\mathcal{S}$  is called the thick envelope of $\mathcal{S}$ and is denoted $\langle \mathcal{S} \rangle$.
One says that $\mathcal{S}$ classically generates $\mathcal{T}$ if $\langle \mathcal{S} \rangle$ is equal to $\mathcal{T}$.
\end{defi}

\begin{defi}
Assume that  $\mathcal{T}$ admits arbitrary direct sums in a given universe. An object $L$ in $\mathcal{T}$ is compact if $\Hom_{\mathcal{T}}(L, \cdot)$ commutes with direct sums. We denote by $\mathcal{T}^{c}$ the full subcategory of $\mathcal{T}$ consisting of compacts objects.
\end{defi}

\begin{defi}
Let $\mathcal{T}$ be a triangulated category admitting arbitrary direct sum in a given universe. The category $\mathcal{T}$ is compactly generated if it is generated by a set of compact objects.
\end{defi}

The following result was proved independently by Ravenel and Neeman, see  \cite{Nee_comp} and \cite{Rav}.

\begin{thm} \label{Rav_Nee}
Assume that $\mathcal{T}$ is compactly generated then a set of objects $\mathcal{S} \subset \mathcal{T}^{c}$ classically generates $\mathcal{T}^{c}$ if and only if it generates $\mathcal{T}$.
\end{thm}

We give an inductive description of the thick envelope of a subset of a triangulated category. For that purpose, we introduce a multiplication on the set of full subcategories of a triangulated category. We follow closely the exposition of \cite{BVdB}. 

\begin{defi}
Let $\mathcal{T}$ be a triangulated category. Let $\mathcal{C}$ and $\mathcal{D}$ be full subcategories of $\mathcal{T}$. Let $\mathcal{C} \star \mathcal{D}$ be the strictly full subcategory of $\mathcal{T}$ whose objects $E$ occur in a triangle of the form 
\begin{equation*}
C \to E \to D \to C \lbrack1 \rbrack
\end{equation*}
where $C \in \mathcal{C}$ and $B \in \mathcal{D}$.
\end{defi}

\begin{prop} \label{asso}
The operation $\star$ is associative.
\end{prop}

Let $\mathcal{S}$ be a set of objects of $\mathcal{T}$. We denote by $\add(\mathcal{S})$ the smallest full subcategory in $\mathcal{T}$ which contains $\mathcal{S}$ and is closed under taking finite direct sums and shifts.\\
We denote by $\smd(\mathcal{S})$ the smallest full subcategory which contains $\mathcal{S}$ and is closed under taking direct summands.

\begin{lemme} \label{rec}
If $\mathcal{C}$ \and $\mathcal{D}$ are closed under finite direct sums, then
$\smd(\smd(\mathcal{C}) \star \mathcal{D})=\smd(\mathcal{C} \star \mathcal{D})$.
\end{lemme}

Denote
\begin{center}
\begin{equation*}
\langle \mathcal{S} \rangle_1=\smd(\add(\mathcal{S}))
\end{equation*}
\begin{equation*}
\langle \mathcal{S} \rangle_k=\smd(\underbrace{\langle \mathcal{S} \rangle_1 \star \ldots \star \langle \mathcal{S} \rangle_1}_{k \; factors})
\end{equation*}
\begin{equation*}
\langle \mathcal{S} \rangle = \bigcup_k \langle \mathcal{S} \rangle_k.
\end{equation*}
\end{center}

Then $\langle \mathcal{S} \rangle$ is the thick envelope of $\mathcal{S}$ (see Definition \ref{epais}).

\subsection{Recollection on algebraic categories}
In this section, we recall some classical facts on algebraic categories, \cite{Keller_deriving}, \cite{Keller_construction}, \cite{Keller_dg}. In this section $R$ is a commutative unital ring.

\begin{defi}
A Frobenius category $\mathcal{E}$ is an exact category (in the sense of Quillen \cite{Quill}) with enough projective and injective objects such that an object is projective if and only if it is injective.
\end{defi}

Let $\sigma$ and $\sigma^{'}$ in $\mathcal{E}$.We denote by $\mathcal{N}(\sigma, \sigma^{'})$ the subgroup of $\Hom_{\mathcal{E}}(\sigma, \sigma^{'})$ formed by the maps that can be factorized through an injective object.
We denote by $\underline{\mathcal{E}}$ the category with the same objects as $\mathcal{E}$ and whose morphisms spaces are the quotients $\Hom_{\mathcal{E}}(\sigma, \sigma^{'})/\mathcal{N}(\sigma, \sigma^{'})$. The category  $\underline{\mathcal{E}}$ is called the stable category of $\mathcal{E}$. A classical result states that $\underline{\mathcal{E}}$ is a triangulated category.

\begin{defi}
One says that an R-linear triangulated category is algebraic if it is triangle equivalent to the stable category of an R-linear Frobenius category.
\end{defi}

\begin{prop}\label{trialg}
A triangulated subcategory of an algebraic triangulated category is algebraic.
\end{prop}

\begin{prop}\label{deralg}
The derived category of a Abelian category is algebraic.
\end{prop}

We have the following theorem from \cite{Keller_construction} which is a consequence of 
\cite[Theorem 4.3] {Keller_deriving}. If $\Lambda$ is a dg category, we denote by $\Der(\Lambda)$ its derived category in the sense of \cite{Keller_deriving} (note that  $\Der(\Lambda)$ is not a dg category).

\begin{thm}\label{equidg}
Let $\mathcal{E}$ be a cocomplete Frobenius category and set $\mathcal{T}=\underline{\mathcal{E}}$. Assume that $\mathcal{T}$ has a compact generator $G$. Then, there is a dg algebra $\Lambda$ and an equivalence of triangulated categories $F: \Der(\Lambda) \to \mathcal{T}$ with $F(\Lambda) \stackrel{\sim}{\to}G$. In particular, we have
\begin{equation*}
\Hn^n(\Lambda) \stackrel{\sim}{\to}\Hom_{\Der(\Lambda)}(\Lambda,\Lambda [n]) \stackrel{\sim}{\to} \Hom_{\mathcal{T}}(G,G[n]),\; n \in \Z. 
\end{equation*}
\end{thm}

\subsection{The case of $\Der_{\qcoh}(\cO_X)$}

 Let $(X, \cO_X)$ be a scheme. We denote by $Qcoh(X)$ the category of quasi-coherent $\cO_X$-modules. Its derived category is denoted by $\Der(Qcoh(X))$. We write $\Der_{\qcoh}(\cO_X)$ for the full triangulated subcategory of $\Der(\cO_X)$ consisting of complexes with quasi-coherent cohomology. 

\begin{thm}[\cite{Nee}] \label{eqqc}
If $X$ is a quasi-compact and separated scheme then the canonical functor $\Der(Qcoh(X)) \to \Der_{\qcoh}(\cO_X)$ is an equivalence.
\end{thm}

\begin{defi}
  Let $(X, \cO_X)$ be a scheme. A perfect complex on $X$ is a complex of $\cO_X$-modules which is locally quasi-isomorphic to a bounded complex of locally free $\cO_X$-modules of finite type. We denote by $\Dper(\cO_X) \subset \Der_{\qcoh}(\cO_X)$ the category of perfect complexes.
\end{defi}

In this paper, we are interested in complex smooth algebraic varieties. We give a few properties of perfect complexes in this setting.  Since $X$ is an algebraic variety, X is a Noetherian topological space. Thus, a perfect complex in $\Der_{\qcoh}(\cO_X)$ is in $\Der^{b}_{\qcoh}(\cO_X)$. Since $\cO_X$ is Noetherian it follows that $\Dper(\cO_X) \subset \Der^{b}_{\coh}(\cO_X)$. Finally since $X$ is smooth, we have $\Der^{b}_{\coh}(\cO_X) \subset \Dper(\cO_X)$. Thus, on a smooth algebraic variety, $\Dper(\cO_X)=\Der^{b}_{\coh}(\cO_X)$.\\

Recall the following theorem from \cite{BVdB}.

\begin{thm}\label{BV}
Assume that $X$ is a quasi-compact and quasi-separated scheme. Then,
\begin{enumerate}[(i)]
\item the compact  objects in $\Der_{\qcoh}(\cO_X)$ are the perfect complexes,
\item $\Der_{\qcoh}(\cO_X)$ is generated by a single perfect complex.
\end{enumerate}
\end{thm}

As a corollary Bondal and Van den Bergh obtain

\begin{thm}\label{dgaff}
Assume that $X$ is a quasi-compact quasi-separated scheme. Then $\Der_{\qcoh}(\cO_X)$ is equivalent to $\Der(\Lambda_0)$ for a suitable dg algebra $\Lambda_0$ with bounded cohomology.
\end{thm}

\subsection{$\hbar$-graduation}\label{sectiongr}

\subsubsection{The case of ringed space} \label{sheafcase}
In this section, X is a topological space and $\mathcal{R}$ is a $\Z [ \hbar]_X$-algebra on $X$ without $\hbar$-torsion.Throughout this text we assume that $\hbar$ is central in $\mathcal{R}$. We set $\cR_0= \cR / \hbar \cR$.
We refer the reader to \cite{KS3} for more details.
\begin{defi}
We denote by $\gr_\hbar: \Der(\cR) \to \Der(\cR_0)$ the left derived functor of the right exact functor $\Mod(\cR) \to \Mod(\cR_0)$ given by $\mathcal{M} \mapsto \mathcal{M} / \hbar \mathcal{M}$. For $\mathcal{M} \in \Der(\cR)$ we call $\gr_\hbar(\mathcal{M})$ the graded module associated to $\mathcal{M}$.We have
\begin{equation*}
\gr_{\hbar} \mathcal{M} \simeq \cR_0 \Lte_\cR \mathcal{M}.
\end{equation*}
\end{defi}

\begin{prop}\label{grhtens}
(i) Let $\cK_1 \in \Der(\op{\cR})$ and $\cK_2 \in \Der(\cR)$. Then,
\begin{equation*}
\gr_\hbar(\cK_1 \Lte_\cR \cK_2) \simeq \gr_\hbar(\cK_1) \Lte_{\cR_0} \gr_\hbar(\cK_2).
\end{equation*}
(ii) Let $\cK_i \in \Der(\cR) \; (i=1,2)$. Then
\begin{equation*}
\gr_\hbar(\fRHom_\cR(\cK_1,\cK_2))\simeq \fRHom_{\cR_0}(\gr_\hbar\cK_1,\gr_\hbar\cK_2).
\end{equation*}
\end{prop}

\begin{prop}
Let $X$ and $Y$ be two topological spaces and $f:X \to Y$ a morphism of topological spaces. The functor $\gr_\hbar : \Der(\Z[\hbar]_X) \to \Der(\Z[\hbar]_Y)$ commutes with the operations of direct images $\dR f_\ast$ and of inverse images. 
\end{prop}

\subsubsection{The case of algebroid stacks}

We write $\C^{\hbar}$ for the ring $\C[[\hbar]]$. In this section $\cA_X$ denotes a $\C^{\hbar}$-algebroid stack without $\hbar$-torsion. As in the previous subsection we refer the reader to \cite{KS3}.

\begin{defi}
Let $\cA_X$ be a $\C^{\hbar}$-algebroid stack without $\hbar$-torsion on a topological space X. One denotes by $\gr_\hbar (\cA_X)$ the $\C$-algebroid associated with the prestack $\mathfrak{S}$ given by:

\begin{flushleft}
$Ob(\mathfrak{S}(U))=Ob(\cA(U))$ for an open set $U$ of $X$,\\
$\Hom_{\mathfrak{S}(U)}(\sigma,\sigma^{'})=\Hom_{\cA}(\sigma,\sigma^{'})/ \hbar \Hom_{\cA}(\sigma,\sigma^{'})$ for $\sigma, \sigma^{'} \in \cA(U)$.
\end{flushleft}
\end{defi}

There is a natural functor $\cA_X \to \gr_\hbar(\cA_X)$ of $\C$-algebroid stacks. This functor induces a functor

\begin{equation*}
\iota_g:\Mod(\gr_\hbar \cA_X) \to \Mod(\cA_X).
\end{equation*}

The functor $\iota_g$ admits a left adjoint functor $\mathcal{M} \to \C \otimes_{\C^{\hbar}} \mathcal{M}$. The functor $\iota_g$ is exact and it induces a functor

\begin{equation*}
\iota_g:\Der(\gr_\hbar \cA_X) \to \Der(\cA_X).
\end{equation*}

On extends the definition of $\gr_\hbar$  by
\begin{equation*}
\gr_\hbar(\mathcal{M}) \simeq \gr_\hbar(\cA_X) \Lte_{\cA_X} \mathcal{M} \simeq \C \Lte_{\C^{\hbar}} \mathcal{M}.
\end{equation*}

The propositions of the preceding subsection concerning sheaves extend to the case of algebroid stacks. Finally we have the following important proposition.

\begin{prop}\label{double_adj}
The functor $\gr_\hbar$ and $\iota_g$ define pairs of adjoint functors $(\gr_\hbar,\iota_g)$ and $(\iota_g,\gr_\hbar[-1])$.
\end{prop}

\begin{proof}
We refer the reader to \cite[proposition 2.3.6]{KS3}. 
\end{proof}

\subsection{Cohomologically Complete Module}
In this subsection, we briefly recall some facts about cohomologically complete modules. We closely follow \cite{KS3} and refer the reader to it for an in depth treatment of the notion of cohomological completeness. 

In this section, $X$ is a topological space and $\mathcal{R}$ is a $\Z [\hbar]_X$-algebra without $\hbar$-torsion.
 We set $\mathcal{R}^{loc}:=\Z[\hbar, \hbar^{-1}]\otimes_{\Z[\hbar]} \mathcal{R}$.

The right orthogonal category $\Der(\mathcal{R}^{loc})^{\bot r}$ to the full subcategory $\Der(\mathcal{R}^{loc})$ of $\Der(\mathcal{R})$ is the full triangulated subcategory consisting of objects $\mathcal{M} \in \Der(\mathcal{R})$ satisfying $\Hom_{\Der(\mathcal{R})}(\mathcal{N},\mathcal{M}) \simeq 0$ for any $\mathcal{N} \in \Der(\mathcal{R}^{loc})$.

\begin{defi}
An object $\mathcal{M} \in  \Der(\mathcal{R})$ is cohomologically complete if its belong to $\Der(\mathcal{R}^{loc})^{\bot r}$. We write $\Der_{\cc}(\cR)$ for  $\Der(\mathcal{R}^{loc})^{\bot r}$.
\end{defi}

 Propositions \ref{ccfact}, \ref{ccgr}, \ref{ccdirect} are proved in \cite{KS3}.

\begin{prop} \label{ccfact}
\begin{enumerate}[(i)]
\item For $\mathcal{M} \in \Der(\mathcal{R})$, the following conditions are equivalent:\\

\begin{enumerate}[(a)]
\item $\mathcal{M}$ is cohomologically complete,\\

\item $\fRHom_{\mathcal{R}}(\mathcal{R}^{loc}, \mathcal{M}) \simeq \fRHom_{\Z[\hbar]}(\Z[\hbar, \hbar^{-1}], \mathcal{M}) \simeq 0,$\\

\item For any $x \in X$, $j=0,1$ and any $i \in \Z$, 
\begin{equation*}
\varinjlim_{x \in U} \Ext^j_{\mathcal{R}}(\mathcal{R}^{loc}, \Hn^i(U,\mathcal{M})) \simeq 0. 
\end{equation*}
Here, $U$ ranges an open neighborhood system of $x$.\\
\end{enumerate}
\item $\fRHom_{\mathcal{R}}(\mathcal{R}^{loc}/ \mathcal{R}, \mathcal{M})$ is cohomologically complete for any $\mathcal{M} \in \Der(\mathcal{R})$.\\
\item For any $\mathcal{M} \in \Der(\mathcal{R})$, there exists a distinguished triangle
\begin{equation*}
\mathcal{M}^{'} \to \mathcal{M} \to \mathcal{M}^{''} \stackrel{+1}{\to}
\end{equation*}
with $\mathcal{M}^{'} \in \Der(\mathcal{R}^{loc})$ and $\mathcal{M}^{''} \in \Der_{\cc}(\cR)$.

\item Conversely, if
\begin{equation*}
\mathcal{M}^{'} \to \mathcal{M} \to \mathcal{M}^{''} \stackrel{+1}{\to}
\end{equation*}
is a distinguished triangle with $\mathcal{M}^{'} \in \Der(\mathcal{R}^{loc})$ and $\mathcal{M}^{''} \in \Der_{\cc}(\cR)$, then $\mathcal{M}^{'} \simeq \fRHom_{\mathcal{R}}(\mathcal{R}^{loc}, \mathcal{M})$ and $\mathcal{M}^{''} \simeq \fRHom_{\mathcal{R}}(\mathcal{R}^{loc}/ \mathcal{R}[-1], \mathcal{M})$.
\end{enumerate}
\end{prop}

\begin{lemme}\label{Homcc}
Assume that $\mathcal{M} \in \Der(\mathcal{R})$ is cohomologically complete. Then $\fRHom_{\cR}(\mathcal{N},\mathcal{M}) \in \Der(\Z_X[\hbar])$ is cohomologically complete for any $\mathcal{N} \in \Der(\cR)$.
\end{lemme}

\begin{prop}\label{ccgr}
Let $\mathcal{M} \in \Der(\mathcal{R})$ be a cohomologically complete object. If $\gr_{\hbar} \mathcal{M} \simeq 0$, then $\mathcal{M} \simeq 0$.
\end{prop}

\begin{cor}\label{isogr}
Let $f:\mathcal{M} \to \mathcal{N}$ be a morphism of $\Der_{\cc}(\cR)$. If $\gr_{\hbar}(f)$ is an isomorphism then $f$ is an isomorphism. 
\end{cor}

\begin{prop}\label{ccdirect}
Let $f:X \to Y$ be a continuous map, and $\mathcal{M} \in \Der(\Z_X[\hbar])$. If $\mathcal{M}$ is cohomologically complete, then so is $\dR f_\ast \mathcal{M}$.
\end{prop}

The following result is implicit in \cite{KS3}. We make it explicit since we use it frequently.

\begin{prop} \label{torsion comp}
Let $\mathcal{M} \in \Der(\cR)$  such that there locally exists $n \in \N$, such that $\hbar^n \mathcal{M}\simeq 0$. Then $\mathcal{M}$ is cohomologically complete.
\end{prop}

\begin{proof}
The question is local. Thus we can assume that $\cA_X$ is a sheaf. The action of $\hbar$ on $\cA_X^{loc}$ is an isomorphism thus the morphism
\begin{equation*}
\hbar \circ: \fRHom(\cA_X^{loc},\mathcal{M}) \to \fRHom(\cA_X^{loc},\mathcal{M})
\end{equation*}
is an isomorphism.
The morphism
\begin{equation*}
\circ \hbar: \fRHom(\cA_X^{loc},\mathcal{M}) \to \fRHom(\cA_X^{loc},\mathcal{M})
\end{equation*}
is locally nilpotent.
Since $\hbar$ is central in $\mathcal{R}$, then $\hbar \circ= \circ \hbar$. Thus, $ \fRHom(\cA_X^{loc},\mathcal{M})=0$.
\end{proof}

\subsection{Modules over formal deformations after \cite{KS3}}

In this subsection, we recall some facts about formal deformation of ringed spaces. We refer the reader to \cite{KS3} for DQ-modules. (Note that they are called twisted deformations in \cite{Yeku1}). We refer to \cite{Stacks_Gerbes}, \cite{categories_and_sheaves} for stacks and algebroid stacks. We denote by $\C^{\hbar}$ the ring $\C[[\hbar]]$.\\

\begin{defi}[\cite{KS3}]
Let $(X, \cO_X)$ be a commutative ringed space on a topological space $X$. Assume that $\cO_X$ is a Noetherian sheaf of $\C$-algebras. A formal deformation algebra $\cA_X$ of $\cO_X$ is a sheaf of $\C^{\hbar}$-algebras such that
\begin{enumerate}[(i)]
\item $\hbar$ is central in $\cA_X$
\item $\cA_X$ has no $\hbar$-torsion
\item $\cA_X$ is $\hbar$-complete
\item $\cA_X / \hbar \cA_X \simeq \cO_X$ as sheaves of $\mathbb{C}$-algebras.
\item There exists a base $\mathfrak{B}$ of open subsets of $X$ such that for any $U \in \mathfrak{B}$ and any coherent $\cO{_X}|_U$-module $\mathcal{F}$, we have $\Hn^n(U,F)=0$ for any $n>0$.
\end{enumerate}
\end{defi}

\begin{Rem}
Clearly, on a complex algebraic variety, condition (iv) of the preceding definition is satisfied.
\end{Rem}

\begin{defi}
A formal deformation algebroid $\cA_X$  on $X$ is a $\C^{\hbar}$-algebroid such that for each open set $U \subset X$ and each $\sigma \in \cA_X(U)$, the $\C^{\hbar}$-algebra $\mathcal{E}nd_{\cA_X}(\sigma)$ is a  formal deformation algebra on U.
\end{defi}

Let $\cA_X$ be a formal deformation algebroid on $X$. We denote by $\Mod(\cA_X)$ the category of functors $\Fct(\cA_X, \sMod(\C^{\hbar}_X))$. The category $\Mod(\cA_X)$ is a Grothendieck category.
For a module $\mathcal{M}$ over an  algebroid $\cA_X$ the local notion of being coherent, locally free etc. still make sense.We denote by $\Der(\cA_X)$ the derived category of $\Mod(\cA_X)$, by $\Der^{b}(\cA_X)$ its bounded derived category and by $\Der^{b}_{\coh}(\cA_X)$ the full triangulated subcategory of $\Der^{b}(\cA_X)$ consisting of objects with coherent cohomologies.

\begin{defi}
We say that an algebroid is trivial if it is equivalent to the algebroid stack associated to a sheaf of rings.
\end{defi}

From now on, we assume that $X$ is a smooth algebraic variety  endowed with the Zarisky topology. There are the following results (see Remark 2.1.17 of \cite{KS3} due to Prof. Joseph Oesterlé)

\begin{prop}
On a smooth algebraic variety $X$, the group $\Hn^2(X,\cO_X^{\times})$ is zero.
\end{prop}

\begin{cor}\label{grgentil}
On a smooth algebraic variety, invertible $\cO_X$-algebroid stacks are trivial.
\end{cor}

By the definition of the functor $\gr_\hbar$, it is clear that $\gr_\hbar \cA_X$ is an $\cO_X$ invertible algebroid (see \cite{KS3} for a definition of invertible) and by Corollary \ref{grgentil} it follows that

\begin{equation}
\gr_\hbar \cA_X \simeq \cO_X. 
\end{equation}

\begin{prop}
The functor $\gr_\hbar$ induces a functor
\begin{equation*}
\gr_\hbar: \Der^{b}_{\coh}(\cA_X) \to \Der^{b}_{\coh}(\cO_X).
\end{equation*}
\end{prop}
We have the following results from \cite{KS3}.

\begin{prop}\label{resloclibre}
Let $d \in \N$. Assume that any coherent $\cO_X$-module locally admits a resolution of length $\leq d$ by free $\cO_X$-modules  of finite rank. Let $\mathcal{M}^{\bullet}$ be a complex of $\cA_X$-modules concentrated in degrees $[a,b]$ and assume that $\Hn^i(\mathcal{M})$ is coherent for all $i$. Then, in a neighborhood of each $x \in X$, there exists a quasi-isomorphism $\mathcal{L}^{\bullet} \to \mathcal{M}^{\bullet}$ where $\mathcal{L}^{\bullet}$ is a complex of free $\cA_X$-modules of finite rank concentrated in degrees $[a-d-1,b]$.
\end{prop}

We have the following sufficient condition which is a corollary of more general results that ensure that under certain conditions, an algebroid stack of formal deformation is trivial (see \cite{Kos}, \cite{Defor_Tsy}, \cite{CalHal}, \cite{Yeku1}).

\begin{prop}\label{deforsheaf}
Let $X$ be a smooth algebraic variety endowed with a deformation algebroid $\cA_X$. If $\Hn^{1}(X,\cO_X)=\Hn^{2}(X,\cO_X)=0$, then $\cA_X$ is equivalent to the algebroid stack associated to a  formal deformation algebra of $\cO_X$.
\end{prop}

\section{Q.C.C modules}

\subsection{Graded quasi-coherent modules and quasi-coherent $\cO_X$-modules}

We start by recalling some results concerning the derived category of quasi-coherent sheaves.

\begin{defi}
Let $\mathcal{M} \in \Der(\cA_X)$.We say that $\mathcal{M}$ is graded quasi-coherent if $\gr_\hbar(\mathcal{M}) \in \Der_{\qcoh}(\cO_X)$. We denote by $\Der_{\gqcoh}(\cA_X)$ the full subcategory of $\Der(\cA_X)$ formed by graded quasi-coherent modules.
\end{defi}

\begin{prop}
The category $\Der_{\gqcoh}(\cA_X)$ is a triangulated subcategory of $\Der(\cA_X)$.
\end{prop}

\begin{proof}
Obvious.
\end{proof}

\subsection{Q.C.C objects}
In this subsection, we introduce the category of $\qcc$-modules.
\begin{defi}
An object $\mathcal{M} \in \Der(\cA_X)$ is $\qcc$ if it is graded quasi-coherent and cohomologically complete. The full subcategory of $\Der(\cA_X)$ formed by $\qcc$-modules is denoted by $\Der_{\qcc}(\cA_X)$.
\end{defi}
Since, $\Der_{\qcc}(\cA_X) = \Der_{\gqcoh}(\cA_X) \cap \Der_{\cc}(\cA_X)$, we have
\begin{prop}\label{qcc_tri}
The category $\Der_{\qcc}(\cA_X)$ is a $\C^{\hbar}$-linear triangulated subcategory of $\Der(\cA_X)$.
\end{prop}

\begin{prop}
If $\mathcal{M} \in \Der^{b}_{\qcc}(\cA_X)$ is such that $\gr_\hbar \mathcal{M} \in \Der^b_{\coh}(\cO_X)$, then $\mathcal{M} \in \Der^b_{\coh}(\cA_X)$.
\end{prop}

\begin{proof}
It is a direct consequence of  \cite[Theorem 1.6.4]{KS3}.
\end{proof}

\begin{prop}
If $\mathcal{M} \in \Der^{b}_{\coh}(\cA_X)$, then $\mathcal{M} \in \Der^b_{\qcc}(\cA_X)$.
\end{prop}

\begin{proof}
It is a direct consequence of  \cite[Theorem 1.6.1]{KS3}.
\end{proof}

We now prove that $\Der_{\qcc}(\cA_X)$ is cocomplete. For that, we first prove that $\Der_{cc}(\cA_X)$ is cocomplete.

\begin{defi}
We denote by $(\cdot)^{\cc}$ the functor 
\begin{equation*}
\fRHom_{\cA_X}((\cA_X^{loc}/\cA_X)[-1],\cdot): \Der(\cA_X) \to \Der(\cA_X). 
\end{equation*}
We call this functor the functor of cohomological completion.
\end{defi}

The following exact sequence
\begin{equation}\label{suiteloc}
0 \to \cA_X \to \cA_X^{loc} \to \cA_X^{loc}/ \cA_X \to 0.
\end{equation}
induces a morphism
\begin{equation*}
\cA_X^{loc}/ \cA_X[-1] \to \cA_X.
\end{equation*}
This morphism yields to a morphism of functor
\begin{equation} \label{morcc}
cc:\id \to (\cdot)^{cc}.
\end{equation}

\begin{prop}\label{cciso}
The morphism of functor 
\begin{equation*}
\gr_\hbar(cc): \gr_\hbar \circ \id \to\gr_\hbar \circ (\cdot)^{cc}
\end{equation*}
is an isomorphism.
\end{prop}

\begin{proof}
We have the following isomorphism
\begin{align*}
\gr_{\hbar}\fRHom_{\cA_X}((\cA_X^{loc}/ \cA_X)[-1],\mathcal{M})&\simeq&\\ 
& \hspace{-4.5mm}\fRHom_{\gr_\hbar \cA_X}(\gr_\hbar (\cA_X^{loc}/ \cA_X)[-1],\gr_\hbar\mathcal{M}).
\end{align*}

Applying the functor $\gr_\hbar$ to (\ref{suiteloc}), we obtain the following distinguished triangle.

\begin{equation*}
\gr_\hbar (\cA_X^{loc}/ \cA_X)[-1]  \to \gr_\hbar \cA_X \to  \gr_\hbar \cA_X^{loc} \stackrel{+1}{\longrightarrow}.
\end{equation*}
Noticing that $ \gr_\hbar \cA_X^{loc} \simeq 0$, we deduce that the map $\gr_\hbar (\cA_X^{loc}/ \cA_X)[-1]  \to \gr_\hbar \cA_X$ is an isomorphism which proves the claim
\end{proof}

\begin{cor} \label{grcomp}
For every $\mathcal{M} \in \Der(\cA_X)$,
\begin{equation*}
\gr_\hbar \mathcal{M}^{\cc} \simeq \gr_\hbar \mathcal{M}.
\end{equation*}
\end{cor}

\begin{defi}
 Let $(\mathcal{M}_i)_{i \in I}$ be a family of objects of $\Der_{\cc}(\cA_X)$. We set 
\begin{equation*}
\obplus \mathcal{M}_i= \left(\bigoplus_{i \in I} \mathcal{M}_i \right)^{\cc}
\end{equation*}
where $\bigoplus$ denote the direct sum in the category $\Der(\cA_X)$. 
\end{defi}

\begin{prop}\label{sumcc}
The category $\Der_{\cc}(\cA_X)$ admits direct sums. The direct sum of the family $(\mathcal{M}_i)_{i \in I}$ is given by $\obplus \mathcal{M}_i$.
\end{prop}

\begin{proof}
 Let $(\mathcal{M}_i)_{i \in I}$ be a family of elements of $\Der_{\cc}(\cA_X)$.
By Proposition \ref{ccfact} $(ii)$, $\obplus \mathcal{M}_i$ is cohomologically complete.

Using the natural transformation (\ref{morcc}) we obtain a morphism

\begin{equation*}
cc:  \bigoplus_{i \in I} \mathcal{M}_i  \to \obplus \mathcal{M}_i.
\end{equation*}

It remains to shows that for all $\mathcal{F} \in \Der_{\cc}(\cA_X)$, $cc$ induces an isomorphism

\begin{equation} \label{iso_int}
 \Hom_{\cA_X}( \obplus \mathcal{M}_i, \mathcal{F})\stackrel{\sim}{\to}  \Hom_{\cA_X}( \bigoplus_{i \in I} \mathcal{M}_i, \mathcal{F}).
\end{equation}

It is enough to prove the isomorphism

\begin{equation}\label{isoplus}
\fRHom_{\cA_X}( \obplus \mathcal{M}_i, \mathcal{F})\to \fRHom_{\cA_X}( \bigoplus_{i \in I} \mathcal{M}_i, \mathcal{F}).
\end{equation}

Since both terms of (\ref{isoplus}) are cohomologically complete by Lemma \ref{Homcc}, it remains to check the isomorphism on the associated graded map. Applying $\gr_\hbar$ to (\ref{isoplus}) and using Lemma \ref{grhtens} $(ii)$ and Proposition \ref{cciso} , we obtain an isomorphism

\begin{equation*}
\fRHom_{\gr_\hbar\cA_X}( \gr_\hbar(\obplus  \mathcal{M}_i), \gr_{\hbar} \mathcal{F})  \stackrel{\sim}{\to}  \fRHom_{\gr_\hbar \cA_X}( \gr_\hbar(\bigoplus_{i \in I}  \mathcal{M}_i),\gr_\hbar \mathcal{F}).
\end{equation*}

which proves the isomorphism (\ref{isoplus}).

Moreover by definition of the direct sum, we have

\begin{equation} \label{final}
 \Hom_{\cA_X}(\bigoplus_{i \in I} \mathcal{M}_i, \mathcal{F}) \stackrel{\sim}{\to} \prod_{i \in I}  \Hom_{\cA_X}( \mathcal{M}_i, \mathcal{F}).
\end{equation}

Composing the isomorphisms (\ref{iso_int}) and (\ref{final}), we obtain the following functorial isomorphism

\begin{equation*}
 \Hom_{\cA_X}( \obplus \mathcal{M}_i, \mathcal{F}) \stackrel{\sim}{\longrightarrow}\prod_{i \in I} \Hom_{\cA_X}( \mathcal{M}_i, \mathcal{F})
\end{equation*}
which prove the proposition.
\end{proof}

\begin{prop}\label{cocoqcc}
The category $\Der_{\qcc}(\cA_X)$ admits direct sums. The direct sum of the family $(\mathcal{M}_i)_{i \in I}$ is given by $\obplus \mathcal{M}_i$.
\end{prop}

\begin{proof}
We know by Proposition \ref{sumcc}, that $\Der_{\cc}(\cA_X)$ admits direct sums and it is given by $\obplus[]$. Let $(\mathcal{M}_i)_{i  \in I} \in \Der_{\qcc}(\cA_X)$. Then, by Corollary \ref{grcomp}, $\gr_\hbar \obplus \mathcal{M}_i =\bigoplus_{i \in I} \gr_\hbar \mathcal{M}_i$. It follows that $\obplus \mathcal{M}_i \in \Der_{\qcc}(\cA_X)$.
\end{proof}

\subsection{Compact objects and generators in $\Der_{\qcc}(\cA_X)$}
In this subsection, we show that $\Der_{\qcc}(\cA_X)$ is generated by a compact generator and we describe its compact objects. We start by proving some additional properties on the functors $\gr_\hbar$ and  $\iota_g$ which are defined in subsection \ref{sectiongr}. 

Concerning $\iota_g$, recall that there is a functor of stacks $\cA_X \to \gr_{\hbar}(\cA_X) \simeq \cO_X$ inducing

\begin{equation*}
\iota_g : \Der(\cO_X) \to \Der(\cA_X)
\end{equation*}

Notice that $\cO_X$ can be endowed with a structure of left $\cO_X$-module and right $\cA_X$-module. When endowed with such structures we denote it by $\cO_{XA}$. The module $\cO_{XA }$ belongs to $\Der^b(\cO_X \otimes_\C \op{\cA_{X}})$. Similarly we have ${_A}\cO_X \in \Der^{b}(\cA_X \otimes_{\C} \op{\cO_X})$. When $\cO_X$ is endowed with its structure of $\cA_X \otimes \op{\cA_X}$-module we denote it by ${_A} \cO_{XA} \in \Der^b(\cA_X \otimes_C \op{\cA_X})$.\\
With these notations, we have
\begin{align*}
\gr_\hbar(\mathcal{M})= \cO_{XA} \Lte_{\cA_X} \mathcal{M}&& \textnormal{and}&&\iota_g(\mathcal{M})= {_A}\cO_X \otimes_{\cO_X} \mathcal{M}.\\
\end{align*}
Hence 
\begin{align*}
\iota_g \circ \gr_\hbar(\mathcal{M})&= {_A}\cO_X \Lte_{\cO_X} \cO_{XA} \Lte_{\cA_X} \mathcal{M}\\
                                                         &\simeq{_A}\cO_{XA} \Lte_{\cA_X} \mathcal{M}.
\end{align*}

\begin{prop}\label{iotafaith}
For every $\mathcal{M} \in \Der(\cO_X)$,
\begin{equation*}
 \iota_g \circ \gr_\hbar \circ \iota_g (\mathcal{M}) \simeq \iota_g(\mathcal{M}) \oplus \iota_g(\mathcal{M})[1]
\end{equation*}
\end{prop}

\begin{proof}
 We have the exact sequence of $\cA_X \otimes \op{\cA_X}$-modules
\begin{equation*}
\xymatrix{0 \ar[r]& \cA_X \ar[r]^-{\hbar}& \cA_X \ar[r]& {_A}\cO_{XA} \ar[r] & 0.}
\end{equation*}
Thus, for every $\mathcal{M} \in \Der(\cA_X)$, we have $\iota_g \circ \gr_\hbar(\mathcal{M}) \simeq  (\cA_X \stackrel{\hbar}{\to} \cA_X) \Lte_{\cA_X} \mathcal{M}$. Hence, for $\mathcal{M} \in \Der(\cO_X)$, $ \iota_g \circ \gr_\hbar \circ \iota_g (\mathcal{M}) \simeq \iota_g(\mathcal{M}) \oplus \iota_g(\mathcal{M})[1]$. 
\end{proof}

\begin{cor}\label{conq}
If $\mathcal{M} \in \Der_{\qcoh}(\cO_X)$, then $\iota_g(\mathcal{M}) \in \Der_{\qcc}(\cA_X)$.
\end{cor}

\begin{proof}
Let $\mathcal{M}$ is in $\Der_{\qcoh}(\cO_X)$ and consider $\gr_\hbar \circ \iota_g(\mathcal{M})$. We compute $\Hn^{i}(\gr_\hbar \circ \iota_g(\mathcal{M}))$.
\begin{align*}
\iota_g(\Hn^{i}(\gr_\hbar \circ \iota_g(\mathcal{M})))& \simeq \Hn^{i}(\iota_g \circ \gr_\hbar \circ \iota_g(\mathcal{M}))\\
& \simeq \Hn^{i}(\iota_g(\mathcal{M}) \oplus \iota_g(\mathcal{M})[1])\\
&\simeq \iota_g(\Hn^{i}(\mathcal{M}) \oplus \Hn^{i+1}(\mathcal{M})).
\end{align*} 
The functor $\iota_g: \Mod(\cO_X) \to \Mod(\cA_X)$ is fully faithful thus
\begin{equation*}
\Hn^{i}(\gr_\hbar \circ \iota_g(\mathcal{M}))\simeq \Hn^{i}(\mathcal{M}) \oplus \Hn^{i+1}(\mathcal{M})
\end{equation*}
thus $\iota_g(\mathcal{M})$ is in $\Der_{\gqcoh}(\cA_X)$ and it is cohomologicaly complete by Proposition \ref{torsion comp}.
\end{proof}

\begin{prop}\label{genqcc}
If $\mathcal{G}$ is a generator of $\Der_{\qcoh}(\cO_X)$, then $\iota_g(\mathcal{G})$ is a generator of $\Der_{\qcc}(\cA_X)$
\end{prop}

\begin{proof}
 By Proposition \ref{conq}, $\iota_g(\mathcal{G})$ is in $\Der_{\qcc}(\cA_X)$.
Let $\mathcal{M} \in \Der_{\qcc}(\cA_X)$ with $\RHom_{\cA_X}(\iota_g(\mathcal{G}),\mathcal{M})=0$. 
By Proposition \ref{double_adj}, we have
\begin{equation*}
 \RHom_{\cA_X}(\iota_g(\mathcal{G}),\mathcal{M}) \simeq  \RHom_{\cO_X}(\mathcal{G}, \gr_\hbar(\mathcal{M})[-1]).
\end{equation*}
Thus, $\RHom_{\cO_X}(\mathcal{G}, \gr_\hbar(\mathcal{M})[-1]) \simeq 0$ and $ \gr_\hbar(\mathcal{M})[-1]$ is in $\Der_{\qcoh}(\cO_X)$ thus $ \gr_\hbar(\mathcal{M})[-1] \simeq 0$.
 Since $\mathcal{M}$ is cohomologically complete, $\mathcal{M} \simeq 0$.
\end{proof}

\begin{lemme}\label{impli}
If $\mathcal{F} \in \Der^{b}_{\coh}(\cA_X)$ satisfies $\cA_X^{loc} \Lte_{\cA_X} \mathcal{F}=0$ then $\mathcal{F}$ is compact in $\Der_{\qcc}(\cA_X)$.
\end{lemme}

\begin{proof}
Let $(\mathcal{M}_i)_{i \in I}$ a family of objects of $\Der_{\qcc}(\cA_X)$. By adjunction in between $(\cA_{X}^{loc}/ \cA_X)[-1] \Lte_{\cA_X} \cdot$ and $\fRHom_{\cA_X}(\cA_{X}^{loc}/ \cA_X)[-1],\cdot)$, we have
\begin{equation*}
\Hom_{\cA_X}(\mathcal{F}, \obplus \mathcal{M}_i) \simeq \Hom_{\cA_X}((\cA_{X}^{loc}/ \cA_X)[-1] \Lte_{\cA_X} \mathcal{F}, \bigoplus_{i\in I} \mathcal{M}_i).  
\end{equation*}

In $\Mod(\cA_X \otimes \op{\cA_{X}})$, we have the exact sequence

\begin{equation*}
0 \to \cA_X \to \cA_{X}^{loc} \to \cA_{X}^{loc} / \cA_X \to 0.
\end{equation*}

Tensoring by $\mathcal{F}$, we obtain the distinguished triangle in $\Der(\cA_X)$

\begin{equation*}
\xymatrix{ \cA_{X}^{loc} / \cA_X[-1] \Lte_{\cA_X} \mathcal{F} \ar[r] &\cA_X \Lte_{\cA_X} \mathcal{F} \ar[r]& \cA_{X}^{loc} \Lte_{\cA_X} \mathcal{F} \ar[r]^-{+1}& \\
}.
\end{equation*}
Since $\cA_{X}^{loc} \Lte_{\cA_X} \mathcal{F}=0$, $\cA_{X}^{loc} / \cA_X[-1] \Lte_{\cA_X} \mathcal{F}$ is isomorphic to $\mathcal{F}$. It follows that

\begin{equation}\label{step1}
\Hom_{\cA_X}(\mathcal{F}, \obplus \mathcal{M}_i) \simeq \Hom_{\cA_X}(\mathcal{F}, \bigoplus_{i\in I} \mathcal{M}_i).
\end{equation}

The module $\mathcal{F}$ belongs to $\Der_{\coh}^{b}(\cA_X)$. Using Proposition \ref{resloclibre}, and the fact that $X$ is a Noetherian topological space, we have
 
\begin{equation}\label{step2}
\Hom_{\cA_X}(\mathcal{F}, \bigoplus_{i\in I} \mathcal{M}_i) \simeq \bigoplus_{i\in I} \Hom_{\cA_X}(\mathcal{F},  \mathcal{M}_i)
\end{equation}
 which together with (\ref{step1}) prove the lemma.
\end{proof}   

\begin{cor}
If $\mathcal{F}$ is compact in $\Der_{\qcoh}(\cO_X)$ then $\iota_g(\mathcal{F})$ is compact in $\Der_{\qcc}(\cA_X)$.
\end{cor}

\begin{cor}\label{compge}
If $\mathcal{G}$ is a compact generator of $\Der_{\qcoh}(\cO_X)$ then $\iota_g(\mathcal{G})$ is a compact generator of $\Der_{\qcc}(\cA_X)$.
\end{cor}

\begin{cor}
The category $\Der_{\qcc}(\cA_X)$ is compactly generated.
\end{cor}

\begin{proof}
By Theorem \ref{BV} due to Bondal and Van den Bergh, $\Der_{\qcoh}(\cO_X)$ has a compact generator. Then, the claim is a direct consequence of Corollary \ref{compge}. 
\end{proof}

\begin{thm}\label{qccompact}
An object $\mathcal{M}$ of $\Der_{\qcc}(\cA_X)$ is compact if and only if $ \mathcal{M} \in \Der^{b}_{\coh}(\cA_X)$ and $\cA_X^{loc} \otimes_{\cA_X} \mathcal{M}=0$.
\end{thm}

\begin{proof}
The condition is sufficient by Lemma \ref{impli}.
Let $\mathcal{G}$ be a compact generator of $\Der_{\qcoh}(\cO_X)$. By Theorem \ref{Rav_Nee}, we know that the set  of compact objects of $\Der_{\qcc}(\cA_X)$ is equivalent to the thick envelope $\langle \iota_g (\mathcal{G}) \rangle$. Let us show that if $\mathcal{F} \in \langle \iota_g (\mathcal{G}) \rangle$ then $\mathcal{F} \in \Der^{b}_{\coh}(\cA_X)$ and $\cA_{X}^{loc} \otimes_{\cA_X} \mathcal{F} =0$. We will proceed by induction.\\

Recall that $\langle \iota_g(\mathcal{G}) \rangle_1=\smd(\add( \iota_g(\mathcal{G}))$ (cf. subsection \ref{envsec}) where $\smd(\add( \iota_g(\mathcal{G}))$ denote the smallest full subcategory of $\Der_{\qcc}(\cA_X)$ containing $\add(\iota_g (\mathcal{G}))$ and closed under taking direct summand. The category $\add(\iota_g (\mathcal{G}))$ is the smallest full subcategory of $\Der_{\qcc}(\cA_X)$ which contains $\iota_g (\mathcal{G})$ and is closed under taking finite direct sums and shifts.

It is clear that if $\mathcal{F} \in \add(\iota_g (\mathcal{G}))$, then $\mathcal{F} \in \Der_{\coh}^b(\cA_X)$ and $\cA_{X}^{loc} \otimes_{\cA_X} \mathcal{F}=0$. If $\mathcal{F} \in \smd(\add(\iota_g(\mathcal{G})))$, then there exist $\mathcal{M} \in \smd(\add(\iota_g(\mathcal{G})))$ such that $\mathcal{F} \oplus \mathcal{M} \in \add(\iota_g(\mathcal{G}))$. Hence, $\cA_X^{loc} \otimes_{\cA_X} \mathcal{F}=0$. We still need to check that for every $i \in \Z$, $\Hn^i(\mathcal{F}) \in \Mod_{\coh}(\cA_X)$. The question is local, so we can choose an open set $U$ of $X$ such that $\cA_X |_U$ is trivial. The sheaf $\Hn^i(\mathcal{F}|_U)$ is a direct summand of the coherent sheaf $\Hn^i((\mathcal{F} \oplus \mathcal{M})|_U)$ and $\cA_X$ is a Noetherian sheaf (see subsection 1.1 of \cite{KS3} for a definition) thus $\Hn^i(\mathcal{F}|_U)$ is coherent.\\

Assume  that for every $k \leq n, \langle \mathcal{G} \rangle_k$ is a subcategory of $\Der_{\coh}^{b}(\cA_X)$ and that for every $\mathcal{F} \in  \langle \mathcal{G} \rangle_k, \; \cA_{X}^{loc} \otimes_{\cA_X} \mathcal{F}=0$. Let $\mathcal{F}$ in $ \langle \mathcal{G} \rangle_{n+1}$. By  Lemma \ref{rec}, $\langle \mathcal{G} \rangle_{n+1}=\smd(\langle \mathcal{G} \rangle_{n} \star \langle \mathcal{G} \rangle_{1})$. Since $\Der^b_{\coh}(\cA_X)$ is a triangulated category, the induction hypothesis implies that  $\langle \mathcal{G} \rangle_{n} \star \langle \mathcal{G} \rangle_{1} \subset  \Der^b_{\coh}(\cA_X)$. It follows that $\mathcal{F}$ is a direct summand of an object of the category $\langle \mathcal{G} \rangle_{n} \star \langle \mathcal{G} \rangle_{1}$. Then, $\mathcal{F} \in \Der^{b}_{\coh}(\cA_X)$ and $\cA_{X}^{loc} \otimes_{\cA_X} \mathcal{F} =0$.
\end{proof}

\subsection{DG Affinity of DQ-modules}

In this subsection we prove that category of $\qcc$ DQ-modules is DG affine. 
\begin{thm}\label{affinity}
Assume $X$ is a smooth complex algebraic variety endowed with a deformation algebroid $\cA_X$. Then, $\Der_{\qcc}(\cA_X)$ is equivalent to $\Der(\Lambda)$ for a suitable dg algebra $\Lambda$ with bounded cohomology.
\end{thm}

\begin{proof}
By Proposition \ref{qcc_tri}, $\Der_{\qcc}(\cA_X)$ is a $\C^{\hbar}$-linear triangulated subcategory of $\Der(\cA_X)$ which is algebraic by Proposition \ref{deralg}. It follows, by Proposition \ref{trialg}, that $\Der_{\qcc}(\cA_X)$ is algebraic. By  Proposition \ref{cocoqcc}, $\Der_{\qcc}(\cA_X)$ is a cocomplete category. Moreover, by Corollary \ref{compge}, $\Der_{\qcc}(\cA_X)$ has a compact generator $\mathcal{G}$. It follows from Theorem \ref{equidg} that $\Der_{\qcc}(\cA_X)$ is equivalent to the derived category of a dg algebra $\Lambda$ such that 
\begin{equation*}
\Hn^n(\Lambda)\simeq \Hom_{\cA_X}(\iota_g(\mathcal{G}),\iota_g(\mathcal{G})[n]), \; n \in \Z.
\end{equation*}
Using the adjunction between $\iota_g$ and $\gr_\hbar[-1]$ and \cite[Lemma 3.3.8]{BVdB}, we get that the cohomology of $\Lambda$ is bounded.
\end{proof}

\section{Q.C.C Sheaves on affine varieties}
We assume that $X$ is a smooth algebraic affine variety. In view of Proposition \ref{deforsheaf}, we assume that $\cA_X$ is sheaf of formal deformations.  We set $A=\Gamma(X,\cA_X)$, $B=\Gamma(X, \cO_X)$ and $a_X: X \to \lbrace pt \rbrace$. As usual we denote by $A_X$ (resp. $B_X$) the constant sheaf with stalk A (resp. B).

\subsection{Preliminary results}

\begin{lemme}
The $A_X$-module $\cA_X$ is flat.
\end{lemme}

\begin{proof}
It is a direct consequence of \cite[Theorem 1.6.5]{KS3}.
\end{proof}

\begin{lemme}\label{Commutation_Rf}
Let $f:X \to Y$ be a morphism of varieties and  let $\mathcal{M} \in \Der(f^{-1}\cA_Y)$ then 
\begin{equation*}
\dR f_{\ast} \mathcal{M}^{cc} \simeq (\dR f_{\ast} \mathcal{M})^{\cc} \; in \; \Der(\cA_Y). 
\end{equation*} 
\end{lemme}

\begin{proof}
It is a direct consequence of \cite[1.5.12]{KS3}
\end{proof}

We recall the following classical result.

\begin{lemme}\label{comcom}
Let $M \in \Der^b(B)$. The canonical morphism 
\begin{equation} \label{morpecom}
 M \to \Rg(X,\cO_X \otimes_{B_X} a_X^{-1} M)
\end{equation}
is an isomorphism.
\end{lemme}

\begin{proof}
If $M$ is concentrated in degree zero, the result follows directly from the equivalence of categories between $Qcoh(\cO_X)$ and $\Mod(B)$. The result extends immediately to the derived category because $\cO_X \otimes_{B_X}\cdot$ is an exact functor and  because $\Rg(X, \cdot)$ is exact on $Qcoh(\cO_X)$ since $X$ is affine.
\end{proof}

If $\cR \to \cR^{'}$ is a morphism of sheaves of rings, we denote by $for_{\cR^{'}}$ the forgetful functor from 
$\Der^b(\cR^{'})$ to $ \Der^b(\cR)$ and $\gr_\hbar ^{\C}$ the functor $\Der(\C^{\hbar}_X) \to \Der(\C_X)$,  $\mathcal{M} \mapsto \C_X \Lte_{\C^{\hbar}_X} \mathcal{M}$.

\begin{prop}\label{oubli}
(i) The diagram
\begin{equation*}
\xymatrix{\Der^b(A_X) \ar[r]^-{for_{A_X}}& \Der^b(\C^{\hbar}_X)\\
\Der^b(\cA_X) \ar[r]^-{for_{\cA_X}} \ar[u]^-{for_{\cA_X}^{A_X}} \ar[d]^-{\gr_\hbar}& \Der^b(\C^{\hbar}_X) \ar @{=}[u] \ar[d]^-{\gr_\hbar^{\C}}\\
\Der^b(\cO_X) \ar[r]^-{for_{\cO_X}} \ar[d]_-{for_{\cO_X}^{B_X}}& \Der^b(\C_X) \ar @{=}[d]\\
\Der^b(B_X) \ar[r]^-{for_{B_X}}& \Der^b(\C_X) 
}
\end{equation*}
 is commutative.\\
(ii) The six forgetful functors $for_{(\cdot)}^{(\cdot)}$ commute with $\Rg(X;\cdot)$.
\end{prop}

\begin{proof}
(i) We start by proving that $for_{\cO_X} \circ \gr_\hbar= \gr_{\hbar}^{\C} \circ for_{\cA_X}$. Let $\mathcal{M} \in \Der^b(\cA_X)$. We have  
\begin{align*}
for_{\cO_X} \circ \gr_\hbar(\mathcal{M})&=\cO_X \Lte_{\cA_X} \mathcal{M}\\
                                                                  &\simeq \C \Lte_{\C^{\hbar}} \cA_X \Lte_{\cA_X} \mathcal{M}\\
                                                                 &\simeq \C \Lte_{\C^{\hbar}} for_{\cA_X}(\mathcal{M})\\
                                                                &\simeq \gr_\hbar^{\C} \circ for_{\cA_X}(\mathcal{M}).
\end{align*}
The other commutation relations are obvious and are left to the reader.\\ 

(ii) Let us prove (ii) for $for_{\cO_X}^{B_X}$. The other cases being similar. The functors $for_{\cO_X}^{B_X}:\Mod(\cO_X) \to \Mod(B_X)$ is exact since $\cO_X$ is flat over $B_X$ and this functor is right adjoints of an exact functor. Thus it preserves injective resolution.
\end{proof}

\begin{prop}
Let $\mathcal{M} \in \Der^b(\cA_X)$. Then, there is an isomorphism in $\Der(B)$
\begin{equation*}
\gr_\hbar \Rg(X,\mathcal{M}) \simeq \Rg(X, \gr_\hbar \mathcal{M}).
\end{equation*}
\end{prop}

\begin{proof}
Notice the claim is true when $\mathcal{M} \in \Der^b(\C^{\hbar}_X)$. Indeed, the functor
$\gr_\hbar: \Der^b(\C^{\hbar})\to\Der^b(\C)$ is given by $\C_X \Lte_{\C^{\hbar}_X}\cdot$ where $\C_X$ is in $\Der^b(\C_X \otimes_{\C_X} \C^{\hbar}_X)$ that is to say in $\Der^b(\C^{\hbar}_X)$. In the category $\Der^b(\C^{\hbar}_X)$, $\C_X$ admits a free resolution given by $\C^{\hbar}_X \stackrel{\hbar}{\to}\C^{\hbar}_X$. Thus we can apply the projection formula and we get the isomorphism
\begin{equation*}
p:\C \Lte_{\C^{\hbar}}\dR a_{X\ast} \mathcal{M} \stackrel{\sim}{\to}\dR a_{X\ast}( \C_X \Lte_{\C^{\hbar}_X} \mathcal{M}).
\end{equation*}
We denote by $for_{\cO_X}$ the forgetful functor from $\Der^b(\cO_X)$ to $\Der^b(\C_X)$. In this proof, we write $\gr_\hbar \mathcal{M}$ for $\cO_X \Lte_{\cA_X} \mathcal{M}$, $\gr_{\hbar}^{B_X} \mathcal{M}$ for $B_X \Lte_{A_X} \mathcal{M}$ and $\gr_{\hbar}^{B} M$ for $B \Lte_A M$ with $M \in \Der^b(A)$ .\\
By restriction, there is a morphism $B_X \Lte_{A_X} \mathcal{M} \to \cO_X \Lte_{\cA_X} \mathcal{M}$. It induces a morphism
\begin{align}\label{piece1}
\Hom_{B_X}( for_{\cO_X}^{B_X}(\gr_\hbar \mathcal{M}), for_{\cO_X}^{B_X}(\gr_\hbar \mathcal{M})) \to& \notag \\ & \hspace{-1cm} \Hom_{B_X}( \gr_{\hbar}^{B_X}  \mathcal{M}, for_{\cO_X}^{B_X}(\gr_\hbar \mathcal{M})).
\end{align}
The coevaluation $a^{-1}_{X} \dR a_{X\ast} \mathcal{M} \to \mathcal{M}$ of the adjunction between $a^{-1}_{X}$ and $\dR a_{X\ast}$ induces a morphism
\begin{align} \label{piece2}
 \Hom_{B_X}(\gr_{\hbar}^{B_X} \mathcal{M}, for_{\cO_X}^{B_X}( \gr_\hbar \mathcal{M})) \to & \notag \\  
& \hspace{-1.7cm} \Hom_{B_X}( \gr_{\hbar}^{B_X}( a^{-1}_X \dR a_{X \ast} \mathcal{M}), for_{\cO_X}^{B_X}( \gr_\hbar \mathcal{M})).
\end{align}
The adjunction between $a_{X}^{-1}$ and $\Rg(X, \cdot)$ gives the following isomorphism
\begin{align} \label{piece3}
\Hom_{B_X}(\gr_{\hbar}^{B_X}( a^{-1}_X \dR a_X \mathcal{M}), for_{\cO_X}^{B_X}(\gr_\hbar \mathcal{M}))& \simeq& \notag \\
 & \hspace{-3.5cm} \Hom_{B}(\gr_{\hbar}^{B} \Rg(X,\mathcal{M}), \Rg(X, for_{\cO_X}^{B_X}(\gr_\hbar \mathcal{M}))). 
\end{align}
and by Proposition \ref{oubli} (ii), there is an isomorphism in $\Der^b(B)$, 
\begin{equation*}
\Rg(X, for_{\cO_X}^{B_X}(\gr_\hbar \mathcal{M})) \simeq \Rg(X, \gr_\hbar \mathcal{M}).
\end{equation*}
Hence, the image of the identity by the composition of the maps (\ref{piece1}), (\ref{piece2}) and (\ref{piece3}) leads to a morphism
\begin{equation*}
pr:B \Lte_{A} \Rg(X,\mathcal{M}) \to \Rg(X, \gr_\hbar \mathcal{M})
\end{equation*}
such that $for_{B}(pr)=p$. Since $for_{B}$ is conservative, $pr$ is an isomorphism.
\end{proof}
\subsection{Q.C.C sheaves on affine varieties}
We define the two functors:
\begin{equation*} 
\Phi: \Der_{\qcc}^b(\cA_X) \to \Der_{\cc}^{b}(A), \; \Phi(\mathcal{M})=\Rg(X,\mathcal{M}) 
\end{equation*}
and
\begin{equation*}
\Psi:\Der_{\cc}^{b}(A) \to \Der_{\qcc}^b(\cA_X), \Psi(M)=(\cA_X \otimes_{A_X} a_X ^{-1}M)^{\cc}.
\end{equation*}

\begin{thm}\label{deforq}
Let X be a smooth affine variety. The functors $\Phi$ and $\Psi$ are equivalences of triangulated categories and are inverses one to each other and the diagram below is quasi-commutative
\begin{equation*}
\xymatrix{ \Der_{\qcc}^b(\cA_X) \ar@<.4ex>[r]^-{\Phi} \ar[d]^-{\gr_\hbar}&  \Der^b_{\cc}(A) \ar@<.4ex>[l]^-{\Psi} \ar[d]^-{\gr_\hbar}\\
\Der_{\qcoh}^b(\cO_X) \ar@<.4ex>[r]^-{\Rg(X,\cdot)} &  \Der^b(B) \ar@<.4ex>[l]^-{\cO_X \otimes_{B_X} \cdot}.}
\end{equation*}
\end{thm}

\begin{proof}
Let $\mathcal{M} \in  \Der_{\qcc}^b(\cA_X)$. By definition, 
\begin{equation*}
\Psi \circ \Phi(\mathcal{M})=\fRHom_{\cA_X}((\cA^{loc}_X / \cA_X) \lbrack -1 \rbrack, \cA_X \Lte_A \Rg(X, \mathcal{M})).
\end{equation*}
By adjunction, we have the morphism of functor $a_X^{-1} \circ \dR a_{X\ast} \to \id$. It follows that we have a morphism $ a_X^{-1} \Rg(X,\mathcal{M}) \to \mathcal{M}$. Tensoring by $\cA_X \Lte_A \cdot$ we get  

\begin{equation}\label{piece}
 \cA_X \otimes_{A_X} a_X^{-1} \Rg(X,\mathcal{M}) \to \cA_X \otimes_{A_X} \mathcal{M}.
\end{equation}
Moreover, 
\begin{equation*}
\begin{split}
\Hom_{\cA_X}( \cA_X \otimes_{A_X} \mathcal{M},\mathcal{M}) &\simeq \Hom_A( \mathcal{M}, \fHom_{\cA_X}(\cA_X,\mathcal{M}))\\
&\simeq \Hom_A( \mathcal{M},\mathcal{M}).
\end{split}
\end{equation*}

Consequently the image of the identity gives a morphism $\cA_X \otimes_{A_X} \mathcal{M} \to \mathcal{M}$. By composing with (\ref{piece}), one obtains a morphism
\begin{equation*}
 \cA_X \otimes_{A_X} a_X^{-1} \Rg(X,\mathcal{M}) \to \mathcal{M}.
\end{equation*}
Applying the functor $( \cdot)^{\cc}$ to the preceding morphism we obtain 
\begin{equation*}
(\cA_X \otimes_{A_X} a_X^{-1} \Rg(X,\mathcal{M}))^{\cc} \to \mathcal{M}^{\cc}.
\end{equation*}
Since $\mathcal{M}$ is cohomologically complete, $\mathcal{M}^{\cc}\simeq \mathcal{M}$. Thus

\begin{equation}\label{premier_iso}
( \cA_X \otimes_{A_X} a_X^{-1} \Rg(X,\mathcal{M}))^{\cc} \to \mathcal{M}.
\end{equation}
Applying $\gr_\hbar$ to the preceding formula, and using the well known equivalence 

\begin{equation*}
\xymatrix{\Der_{\qcoh}^b(\cO_X) \ar@<.4ex>[r]^-{\Rg(X,\cdot)} &  \Der^b(B) \ar@<.4ex>[l]^-{\cO_X \otimes_{B_X} -}},
\end{equation*}

 we obtain the isomorphism

\begin{equation*}
\cO_X \otimes_{B_X} \Rg(X, \gr_\hbar \mathcal{M}) \stackrel{\sim}{\longrightarrow} \gr_\hbar \mathcal{M}.
\end{equation*}

Since $( \cA_X \otimes_{A_X} a_X^{-1} \Rg(X,\mathcal{M}))^{\cc}$ and $\mathcal{M}$ are cohomologically complete modules, it follows that (\ref{premier_iso}) is an isomorphism.\\

Let $M \in \Der^b_{\cc}(A)$. By definition,
\begin{equation*}
\Phi \circ \Psi(M)=\Rg(X,(\cA_X \otimes_{A_X} a_X^{-1}M)^{\cc}).
\end{equation*}
and using Lemma \ref{Commutation_Rf} we get that
\begin{equation*}
 \Phi \circ \Psi(M) \simeq (\Rg(X,\cA_X \otimes_{A_X} M))^{\cc}.
\end{equation*}

We have a morphism

\begin{equation*}
\dR a_{X\ast} \cA_X \otimes_{A_X} M \to \dR a_{X\ast}( \cA_X \otimes_{A_X} a_X^{-1} M).
\end{equation*}

Since $X$ is affine we obtain $\dR a_{X\ast} \cA_X \simeq A$ thus

\begin{equation*}
M \to \Rg(X, \cA_X \otimes_{A_X} M).
\end{equation*}

We have a map

\begin{equation}\label{equi2}
M^{\cc} \to  (\Rg(X, \cA_X \otimes_{A_X} M))^{\cc}.   
\end{equation}

Applying the functor $\gr_{\hbar}$, we obtain

\begin{equation}\label{gr_proj}
\gr_{\hbar}M \to \Rg(X, \cO_X \otimes_{B_X}  \gr_{\hbar} M).
\end{equation}

Using Lemma \ref{comcom}, we deduce that the map (\ref{gr_proj}) is an isomorphism. It follows by Corollary \ref{isogr} that the morphism (\ref{equi2}) is an isomorphism. This proves the announced equivalence.
\end{proof}

\bibliographystyle{plain}
\bibliography{biblio}

\end{document}